\documentclass[preprint,12pt]{elsarticle}
\usepackage{natbib}
\usepackage{graphicx}
\usepackage{mathtools}
\usepackage{amsthm}
\usepackage{amsmath}
\usepackage{amssymb}
\usepackage{bm}
\usepackage{bbm}
\usepackage{upgreek}
\usepackage[colorlinks=true,linkcolor=blue,citecolor=blue]{hyperref}
\usepackage{enumitem}
\usepackage{multirow}
\usepackage{mathrsfs}
\usepackage{tabularx}

\usepackage{epstopdf}
\usepackage{comment}
\usepackage{color}
\usepackage{graphicx}

\newtheorem{theorem}{Theorem}

\newtheorem{corollary}{Corollary}
\newtheorem{lemma}{Lemma}

\newcommand{\indep}{\raisebox{0.05em}{\rotatebox[origin=c]{90}{$\models$}}}

\newcommand{\X}{{\bf X}}
\newcommand{\Z}{{\bf Z}}
\newcommand{\B}{{\bf B}}
\newcommand{\spc}{{\mathcal S}_{Y|\X}}
\def\R{\mathbb{R}}

\newcommand{\A}{{\bf A}}
\newcommand{\K}{{\bf K}}
\newcommand{\F}{{\bf F}}

\newcommand{\bH}{{\bf H}}
\newcommand{\bP}{{\bf P}}

\newcommand{\h}{{\bf h}}

\newcommand{\s}{{\bf s}}

\newcommand{\tX}{\tilde{\bf X}}
\newcommand{\tbeta}{\tilde{\bm{\beta}}}
\newcommand{\tSigma}{\tilde{\bm{\Sigma}}}

\newcommand{\bdelta}{\bm{\delta}}
\newcommand{\bbeta}{\bm{\beta}}
\newcommand{\BXi}{\bm{\Xi}}
\newcommand{\var}{\text{Var}}
\newcommand{\BSigma}{\bm{\Sigma}}
\newcommand{\BLambda}{\bm{\Lambda}}
\newcommand{\tlambda}{\tilde{\lambda}}
\newcommand{\cov}{\text{Cov}}

\newcommand{\bvarphi}{\bm{\varphi}}

\newcommand{\mI}{\mathbb{I}}

\newcommand{\vect}{\text{Vec}}
\newcommand{\BOmega}{\bm{\Omega}}
\newcommand{\Psibf}{\bm{\Psi}}

\newcommand{\Ib}{{\bf I}}

\newcommand{\Qb}{{\bf Q}}
\newcommand{\Jb}{{\bf J}}
\newcommand{\Kb}{{\bf K}}

\def\eop
{\hfill $\Box$
}

\journal{Statistica Sinica}
\begin{document}
\begin{frontmatter}
\title{On sufficient dimension reduction via principal asymmetric least squares}
\author{Abdul-Nasah Soale}
\author{Yuexiao Dong\corref{cor1}}
\address{Department of Statistical Science, Temple University,
Philadelphia, PA, US, 19122}
 \ead{ydong@temple.edu}
 \cortext[cor1]{Corresponding author.}

\begin{abstract}
In this paper, we introduce principal asymmetric least squares (PALS) as a unified framework for linear and nonlinear sufficient dimension reduction. 
Classical methods such as sliced inverse regression (Li, 1991) and principal support vector machines (Li, Artemiou and Li, 2011) may not perform well in the presence of heteroscedasticity, while our proposal addresses this limitation by synthesizing different expectile levels. Through extensive numerical studies, we demonstrate the superior performance of PALS in terms of both computation time and estimation accuracy. 
For the asymptotic analysis of PALS for linear sufficient dimension reduction, we develop new tools to 
compute the derivative of an expectation of a non-Lipschitz function. 
% We also apply PALS to detect influential observations of the Boston housing data in a model-free fashion. 

%, a useful method for linear sufficient dimension reduction with extension to nonlinear scenarios.  Most existing statistical regression methods do not address heteroscedasticity, which is very common in real data.  Quantile regression is among the few methods in the literature that attempts to deal with the problem. However, quantile regression suffers some known disadvantages.  To escape these disadvantages, we propose PALS which is based on expectile regression.  Expectile regression is a close surrogate to quantile regression but its estimates are more efficient, and it is computationally faster. We demonstrate the competitive performance of the new estimates through extensive numerical studies and real data analysis.
\end{abstract}
\begin{keyword}
distance correlation \sep expectile regression  \sep  heteroscedasticity \sep  nonlinear dimension reduction.
\end{keyword}
\end{frontmatter}

\section{Introduction}

For univariate response $Y$ and multivariate predictor $\X\in \R^p$, sufficient dimension reduction (Li, 1991; Cook, 1998)
aims to find $\B  \in \R^{p \times d}$ such that
\begin{equation}
   Y \indep \X | \B^\top  \X,
    \label{SDR}
 \end{equation}
where ``\indep'' means statistical independence. Under (\ref{SDR}), the conditional distribution of $Y$ given $\X$
is the same as the conditional distribution of $Y$ given $\B^\top  \X$. If $\B$ satisfies (\ref{SDR}), the column space of $\B$ is called a dimension reduction
space. Under very general conditions as discussed in Yin, Li and Cook (2008), the intersection of all dimension reduction spaces is also a dimension reduction space. We refer to this minimum dimension reduction space as the central space for the regression between $Y$ and $\X$, and we denote it as $\spc$. The dimensionality of the central space is known as the structural dimension.

%In the era of big data,  there has been a growing demand for methodology that can extract information from the data effectively.
%Since its inception in Li (1991) and Cook (1998), sufficient dimension reduction (SDR) has seen applications in diverse fields such as bioinformatics (Li and Li, 2004; Liu, Chiaromonte and Li, 2017) and machine learning (Fukumizu, Bach and Jordan, 2004; Suzuki and Sugiyama, 2013; Chen et al., 2018).
%SDR is a supervised dimensionality reduction paradigm that fully retains the regression information between the response and the predictor through projection of the predictor onto a subspace.
%More specifically, for univariate response $Y$ and multivariate predictor $\bm{X}\in \mathbb{R}^p$, SDR 
% where "\indep" means statistical independence. We refer to the corresponding column space as the central space for the regression of $Y$ on $\bm{X}$, and denote it as $\spc$. Under very general conditions (Yin, Li and Cook, 2008), the central space is well-defined and is unique. The dimensionality of the central space is known as the structural dimension.

% Suzuki, T. and Sugiyama, M. (2013) Suﬃcient dimension reduction via squared-loss mutual information estimation. Neural Computation, 25, 725–758.
% \cite{cook1996graphics}

Moment-based methods such as sliced inverse regression (SIR) (Li, 1991), sliced average variance estimation (SAVE) (Cook and Weisberg, 1991), and directional regression (Li and Wang, 2007) are among the most popular sufficient dimension reduction methods. These moment-based methods are easy to implement in practice, and their extensions include sparse sufficient dimension reduction (Li, 2007), dimension reduction with matrix-valued predictors (Li, Kim and Altman, 2010), and dimension reduction for functional data (Li and Song, 2017). For an excellent review, please refer to Li (2018). More recently, Li, Artemiou and Li (2011) proposed the principal support vector machine (PSVM), which applies a modified support vector machine to find the optimal separating hyperplanes of the discretized response. It is shown that the normal vector of the separating hyperplanes can be used to recover $\spc$.
%   By employing support vector machine and focusing on separating hyperplanes rather than slice means, PSVM improves the accuracy of popular inverse regression estimators. 
   Extensions of PSVM include $\ell q$  PSVM (Artemiou and Dong, 2016), principal logistic regression (Shin and Artemiou, 2017), and weighted PSVM (Shin et al., 2017).

A well-known limitation of the moment-based sufficient dimension reduction methods as well as PSVM is that they may not perform well in the presence of heteroscedastic error. In this paper, we propose to replace the hinge loss in PSVM with the asymmetric least squares loss, and we refer to the new proposal as principal asymmetric least squares (PALS). By synthesizing different expectile levels, PALS can improve the performance of PSVM when the error is heteroscedastic. We implement the sample level estimation of PALS through quadratic programming,  provide the asymptotic normality of the sample PALS estimator, and extend PALS for nonlinear sufficient dimension reduction.   Wang, Shin, and Wu (2018) proposed principal quantile regression (PQR), where quantile regression was used instead of the expectile regression in our proposal. We note that the check loss from the PQR objective function is not smooth, while PALS utilizes a smooth objective function. As a result, PALS leads to more accurate estimation with much improved computational speed compared to PQR.

The rest of the paper is organized as follows. The population level development and the sample level estimation of PALS are studied in Section 2 and Section 3, respectively.  Extensions to nonlinear sufficient dimension reduction are examined in Section 4. 
Extensive simulation studies are reported in Section 5 and we provide a real data analysis in Section 6. Section 7 concludes the paper with 
some discussions.  All the proofs are relegated to the Appendix.

\section{Population level development}

Let $\bm{\mu}=E(\X)$ and  $\BSigma=\var(\X)$. The population $\uptau$-th level objective function of PALS is 
\begin{equation}
L_\uptau(\alpha,\bm{\beta}) =  \bm{\beta}^\top   \bm{\bm{\Sigma}} \bm{\beta} + \lambda E\left[ \rho_\uptau\left\{Y-\alpha-\bm{\beta}^\top  (\X-\bm{\mu} )\right\} \right].
\label{obj_pop}
\end{equation}
Here $\uptau \in (0,1)$ denotes the expectile level, $\lambda>0$  is a tuning parameter, and $\rho_\uptau$ is the asymmetric least squares loss function (Newey and Powell, 1987) defined as follows
\begin{align}
\rho_\uptau(c) = \begin{cases}
(1-\uptau)c^2 & \mbox{if } c\leq 0,\\
\uptau c^2 & \mbox{if } c>0.
\end{cases}
\label{als_loss}
\end{align}

The asymmetric least squares loss was originally designed to recover the regression expectiles, which is known be closely related to regression quantiles (Abdous and Remillard, 1995). 
The objective function (\ref{obj_pop}) is linked to sufficient dimension reduction through the next result. 

\begin{theorem}\label{thm_linearity}
Suppose $E(\X|\B^\top  \X)$ is linear in $\B^\top  \X$, where $\B\in \R^{p\times d}$ is a basis of $\spc$. Let
\begin{align*}
(\alpha_{0,\uptau},\bm{\beta}_{0,\uptau}) = \underset{\alpha\in \R, \ \bm{\beta}\in \R^p}{\operatorname{argmin}}\  L_\uptau(\alpha, \bm{\beta}).
\end{align*}
 Then
$\bm{\beta}_{0,\uptau}\in \spc$.
\end{theorem}

\noindent The assumption about $E(\X|\B^\top  \X)$ is known as the linear conditional mean condition, and is common in the sufficient dimension reduction literature. 
As a result of Theorem \ref{thm_linearity}, we have

\begin{corollary}\label{co}
Suppose $E(\X|\B^\top  \X)$ is linear in $\B^\top  \X$, where $\B\in \R^{p\times d}$ is a basis of $\spc$. Let 
$0<\uptau_1<\ldots<\uptau_K<1$ and
\begin{align*}
(\alpha_{0,\uptau_k},\bm{\beta}_{0,\uptau_k}) = \underset{\alpha\in \R, \ \bm{\beta}\in \R^p}{\operatorname{argmin}}\  L_{\uptau_k}(\alpha, \bm{\beta})\mbox{ for }k=1,\ldots,K.
\end{align*}
Then $\mathrm{span}(\BLambda)\subseteq\spc$, where $\BLambda = \sum_{k=1}^K \bm{\beta}_{0,\uptau_k}\bm{\beta}_{0,\uptau_k}^\top$.
\end{corollary}

\noindent Here $\mathrm{span}$ denotes the column space. Corollary \ref{co} suggests that we can recover the central space by optimizing the PALS objective function (\ref{obj_pop}) at multiple expectile levels.

\section{Sample level estimation}
Given an i.i.d sample $\{(\X_i,Y_i):i=1,\ldots,n\}$, the sample version of (\ref{obj_pop}) becomes
\begin{align}
\hat{L}_\uptau(\alpha,\bm{\beta}) =  \bm{\beta}^\top   \hat{\bm{\bm{\Sigma}}} \bm{\beta} + \cfrac{\lambda}{n}\sum_{i=1}^n  \rho_\uptau\left\{Y_i-\alpha-\bm{\beta}^\top  \left(\X_i-\bar{\X} \right)\right\},
\label{obj_sam1}
\end{align}
where $\bar{\X}=n^{-1}\sum_{i=1}^n \X_i$ and $\hat{\BSigma}=n^{-1}\sum_{i=1}^n (\X_i-\bar{\X})(\X_i-\bar{\X})^\top$. 
%are respectively, the sample mean and sample covariance of $\bm{X}$. To simplify (\ref{obj_sam1}), we let $\bm{Z}_i=\hat{\bm{\bm{\Sigma}}}^{-1/2}(\bm{X}_i-\bar{\bm{X}})$ and $\bm{\theta}=\hat{\bm{\bm{\Sigma}}}^{1/2}\bm{\beta}$ which reduces it to
Denote $\tilde{\lambda}=n^{-1}\lambda$, $\Z_i=\hat{\BSigma}^{-1/2}(\X_i-\bar{\X})$ and $\bm{\theta}=\hat{\BSigma}^{1/2}\bm{\beta}$. (\ref{obj_sam1})  reduces to 
\begin{align}
\tilde{L}_\uptau(\alpha,\bm{\theta}) =  \bm{\theta}^\top   \bm{\theta} + \tilde{\lambda}\sum_{i=1}^n  \rho_\uptau(Y_i-\alpha-\bm{\theta}^\top   \Z_i).
\label{obj_sam2}
\end{align}
Let $c_{+}=max(0,c)$. We now introduce 
$$\xi_{i+} = (Y_i - \alpha- \bm{\theta}^\top   \Z_i)_{+} \mbox{ and }\xi_{i-} = ( \alpha+ \bm{\theta}^\top   \Z_i-Y_i)_{+}.$$
From the definition of $\rho_{\uptau}$ in (\ref{als_loss}), 
(\ref{obj_sam2}) leads to the following primal optimization problem
\begin{align}
(\hat\alpha_{0,\uptau},\hat{\bm{\theta}}_{0,\uptau})=\underset{\alpha\in \R, \bm{\theta}\in \R^p}{\operatorname{argmin}} \hspace{0.11in} \bm{\theta}^\top  \bm{\theta} + \tilde{\lambda}\uptau \sum_{i=1}^n \xi_{i+}^2 + \tilde{\lambda}(1-\uptau) \sum_{i=1}^n \xi_{i-}^2
\label{opt_primal}
\end{align}
subject to $\xi_{i+} \geq 0$, $\xi_{i-}\geq 0$,
$\xi_{i+} \geq Y_i - \alpha- \bm{\theta}^\top   \Z_i$,
and $\xi_{i-} \geq \alpha+ \bm{\theta}^\top   \Z_i-Y_i$.

\begin{theorem}\label{thm_duality}
Let $\mathbb{Y}=(Y_1,\ldots,Y_n)^\top  $ and $\mathbb{Z}=(\Z_1^\top  ,\ldots,\Z_n^\top  )^\top  $. 
The dual optimization problem of (\ref{opt_primal})
is
\begin{align}
\begin{split}
 (\hat{\bm{a}}_{0,\uptau},\hat{\bm{\eta}}_{0,\uptau})=\underset{\bm{a}\in \R^n, \bm{\eta}\in \R^n}{\operatorname{argmax}}  \hspace{0.1in}
&(\bm{a}-\bm{\eta})^\top   \mathbb{Y}-\cfrac{1}{4}(\bm{a}-\bm{\eta})^\top  \mathbb{Z}\mathbb{Z}^\top  (\bm{a}-\bm{\eta})
-\cfrac{1}{4\tlambda\uptau}\bm{a}^\top  \bm{a} \\
& \hspace{0.3in} 
-\cfrac{1}{4\tlambda(1-\uptau)}\bm{\eta}^\top  \bm{\eta}
\end{split}
\label{opt_dual}
\end{align}
subject to $\bm{a} \geq \bm{0}_n$,
$\bm{\eta} \geq \bm{0}_n$, and
$(\bm{a}-\bm{\eta})^\top  \bm{1}_n = 0$. Furthermore, we have
\begin{align}
\hat{\bm{\theta}}_{0,\uptau}=\cfrac{1}{2} \mathbb{Z}^\top (\hat{\bm{a}}_{0,\uptau}-\hat{\bm{\eta}}_{0,\uptau}).
\label{pd}
\end{align}
\end{theorem}

Consider expectile levels $0<\uptau_1<\ldots<\uptau_K<1$. For a given $\uptau_k$, the dual problem  (\ref{opt_dual}) can be solved through standard quadratic programming to get $\hat{\bm{a}}_{0,\uptau_k}$ and $\hat{\bm{\eta}}_{0,\uptau_k}$. After computing $\hat{\bm{\theta}}_{0,\uptau_k}$ from (\ref{pd}), we 
get the minimizer of 
$\hat{L}_{\uptau_k}$ in (\ref{obj_sam1}) as $\hat{\bm{\beta}}_{0,\uptau_k}=\hat{\BSigma}^{-1/2}\hat{\bm{\theta}}_{0,\uptau_k}$. Based on Corollary \ref{co}, we get the estimator of $\BLambda=\sum_{k=1}^K {\bm{\beta}}_{0,\uptau_k}{\bm{\beta}}_{0,\uptau_k}^\top$ as $\hat{\BLambda} = \sum_{k=1}^K \hat{\bm{\beta}}_{0,\uptau_k}\hat{\bm{\beta}}_{0,\uptau_k}^\top$. Recall that $d$ denotes the
structural dimension of $\spc$. The eigenvectors corresponding to the $d$ largest eigenvalues of $\hat{\BLambda}$ then consist the final PALS estimator to recover the central space. 

We conclude this section with the asymptotic normality of $\vect({\hat{\BLambda}})$, where $\vect$ means vectorization. The details are provided in the Appendix. In order to compute the derivative of an expectation of a non-Lipschitz function, we extend the theoretical development of PSVM (Li, Artemiou and Li, 2011). This extension is necessary as PSVM deals with discretized response while PALS applies to the continuous response without discretization. 

\begin{theorem}\label{normal}
Suppose the regularity conditions in
Theorem \ref{d} and Theorem \ref{h} from the Appendix are satisfied. Then we have
\begin{align*}
    \sqrt{n} \left\{\vect({\BLambda})-\vect({\hat{\BLambda}})\right\}  \overset{D}{\longrightarrow} N (\bm{0},\BOmega)
\end{align*}
as $n\rightarrow \infty$, where ``$\overset{D}{\longrightarrow}$'' means converge in distribution and $\BOmega$ is specified in the Appendix.
\end{theorem}

\section{Nonlinear sufficient dimension reduction}

Suppose $\bvarphi:\R^p \mapsto \R^d$ with $d<q$ are nonlinear functions satisfying
\begin{align}
Y \indep \X | \bvarphi (\X) \label{model:nonlinear},
\end{align}
where $\bvarphi (\X)=\{\varphi_1(\X),\ldots,\varphi_d(\X)\}$.
Then the conditional distribution of $Y$ given $\X$
is the same as the conditional distribution of $Y$ given $\bvarphi (\X)$, and identifying $\bvarphi (\X)$ is known as nonlinear sufficient dimension reduction.
%Following the development of Li, Artemiou and Li (2011), we discuss PALS for nonlinear sufficient dimension reduction in this section. 
Let $\mathcal{H}$ be a reproducing kernel Hilbert space of the
functions of $\X$ with inner product $\langle \cdot, \cdot
\rangle_{\mathcal{H}}$. Let $\Sigma: \mathcal{H} \mapsto
\mathcal{H}$ be the covariance operator such that
 $\langle f_1,
\Sigma f_2 \rangle_{\mathcal{H}}=\cov\{ f_1 (\X), f_2 (\X) \}$ for
any $f_1, f_2\in \mathcal{H}$.  Consider objective function 
\begin{equation}
\Pi_\tau(\alpha,\varphi) =  \langle \varphi,\Sigma \varphi \rangle_\mathcal{H}  + \lambda E[ \rho_\tau\{Y-\alpha-\varphi(\X)\}].
\label{obj_nonlinear1}
\end{equation}
Compared with (\ref{obj_pop}), we see that $\Pi_\tau(\alpha,\varphi)$ is a generalization of $L_\uptau(\alpha,\bm{\beta})$ with the matrix $\BSigma$ replaced by the operator $\Sigma$, the linear function $\bm{\beta}^\top  \X$ replaced by the nonlinear function $\varphi(\X)$, and the inner product in $\R^p$ replaced by the inner product in ${\mathcal{H}}$. Let $(\alpha_{0,\uptau},\varphi_{0,\uptau})$ be the minimizer of $\Pi_\tau(\alpha,\varphi)$ over $\alpha\in\R$ and $\varphi\in
\mathcal{H}$. Under proper conditions, it can be shown that $\varphi_{0,\uptau}$ is a function of $\varphi_1,\ldots,\varphi_d$. See, for example, Theorem 2 of 
Li, Artemiou and Li (2011).

Based on an i.i.d. sample $\{(\X_i, Y_i): i=1,\ldots,n\}$, we now describe the implementation of nonlinear dimension reduction through PALS. Suppose  $\mathcal{H}$ can be spanned by $\{\psi_1,\ldots, \psi_m\}$. Then any function $\varphi\in \mathcal{H}$ becomes
$\varphi(\X)=\bm{\gamma}^\top   \bm{\psi}(\X)$,
where $\bm{\gamma}\in\R^m$ and 
$\bm{\psi}(\X) = \{\psi_1(\X),\ldots,\psi_m(\X)\}^\top$. The sample version of (\ref{obj_nonlinear1}) thus becomes

\begin{equation}
\hat\Pi_\tau(\alpha,\bm{\gamma}) 
= \frac{1}{n} \bm{\gamma}^\top \Psibf^\top \Psibf \bm{\gamma} + 
\frac{\lambda}{n} \sum_{i=1}^{n}  \rho_\tau\{Y_i-\alpha-\bm{\gamma}^\top   \bm{\psi}(\X_i)\},
\label{obj_nonlinear2}
\end{equation}
where $ \Psibf\in\R^{n\times m}$ and the $i$th row of $\Psibf$ is $\bm{\psi}^\top(\X_i)$.
$\hat\Pi_\tau(\alpha,\bm{\gamma}) $ has the same form as 
$\hat{L}_\uptau(\alpha,\bm{\beta})$  in (\ref{obj_sam1}), and can be minimized in a similar fashion. Denote the minimizer of $\hat\Pi_\tau(\alpha,\bm{\gamma})$
as  $(\hat\alpha_{0,\uptau},\hat{\bm{\gamma}}_{0,\uptau})$. 
We then estimate $\varphi_{0,\uptau}(\X)$ by $\hat{\varphi}_{0,\uptau}(\X)=\hat{\bm{\gamma}}_{0,\uptau}^\top   \bm{\psi}(\X)$.
To synthesize multiple expectile levels, consider expectile levels $0<\uptau_1<\ldots<\uptau_K<1$. For a given $\uptau_k$, we get $\hat{\bm{\gamma}}_{0,\uptau_k}$ from minimizing  $\hat\Pi_{\tau_k}(\alpha,\bm{\gamma})$. 
Denote $\hat {\bm\Gamma}=\sum_{k=1}^K \hat{\bm{\gamma}}_{0,\uptau_k} \hat{\bm{\gamma}}_{0,\uptau_k}^\top$ with $d$ leading eigenvectors as 
$\hat{\bm\nu}_1,\ldots,\hat{\bm\nu}_d$. The final estimator of $\bvarphi (\X)$ in (\ref{model:nonlinear}) is $\{\hat{\bm\nu}_1^\top   \bm{\psi}(\X),\ldots,\hat{\bm\nu}_d^\top   \bm{\psi}(\X)\}$.

It remains to choose a proper basis  $\{\psi_1,\ldots, \psi_m\}$ for $\mathcal{H}$.
Define kernel matrix $\Kb_n\in\R^{n\times n}$, with the element in the $i$th row and $j$th column as
\begin{equation}
\kappa(\X_i,\X_j)=exp(-r\|\X_i-\X_j \|^2).
\label{rbf}
\end{equation}
Here $r$ is a tuning parameter and $\|\cdot\|$ denotes the Euclidean norm.
Define $\Qb_n=\Ib_n-\Jb_n/n$, where $\Ib_n$ is the ${n\times n}$ identity matrix and $\Jb_n$ is the ${n\times n}$ matrix whose entries are $1$. 
%Let $\wbf_g$ be the eigenvector corresponding to $\lambda_g$, the $g$th largest eigenvalue of $\Qb_n \Kb_n\Qb_n$ for $g=1,\ldots, n$.
%From Proposition 2 in Li, Artemiou and Li (2011), we know  $\Psibf$ becomes $(\wbf_1,\ldots,\wbf_G)$.
For 
$j=1,\ldots,m$, let $\bm{w}_j=(w_{j1},\ldots,w_{jn})^\top$ be the eigenvector corresponding to the $j$-th largest eigenvalue of  $\Qb_n \Kb_n\Qb_n$. Then $\psi_j(\X)=\sum_{\ell=1}^n\kappa(\X,\X_\ell) w_{j\ell}$ following Proposition 2 of Li, Artemiou and Li (2011). In our simulations, we choose $m=n/2$, and use  the sample version of $E^{-1/2}(\|\X-\X'\|)$ for $r$ in (\ref{rbf}), where $\X$ and $\X'$ are independent $N(\bf{0},\Ib_p)$.

\section{Simulation studies}

\subsection{Linear sufficient dimension reduction}

We evaluate the performance of PALS for linear sufficient dimension reduction in this section. The following models are considered: 
\begin{align*}
\textrm{I: } & Y = \frac{X_1}{0.5 + (X_2 + 1.5)^2} + \varepsilon;\\
\textrm{II: } & Y = 3\sin\{0.25(X_1+X_2)\} + 3\sin\{0.25(X_3+X_4)\} + \varepsilon;\\
\textrm{III: } & Y =  X_1 +0.5\big(e^{0.15X_2}\big)\varepsilon,
\end{align*}
where $\varepsilon\sim N(0,1)$ and $\varepsilon$ is independent of $\X=(X_1,\ldots,X_p)^\top$. The distribution of $\X$ will be specified later.
Let $\bbeta_1=(1,0,\ldots,0)^\top$, $\bbeta_2=(0,1,0,\ldots,0)^\top$, $\bbeta_3 = (1,1,0,\ldots,0)^\top $, and $\bbeta_4 = (0,0,1,1,0,\ldots,0)^\top$. Denote 
$\B$ as the basis of the central space $\spc$. Then $\B=(\bbeta_1,\bbeta_2)$ for models I and III, while  $\B=(\bbeta_3,\bbeta_4)$ for model II. 

We compare PALS with five existing methods in the literature: SIR, SAVE, directional regression (DR), PSVM, and PQR. The number of slices for SIR is set as $10$, and we use $4$ slices for SAVE and DR. Note that SIR is generally not sensitive to the choice of slice numbers, while SAVE and DR work better with fewer slices. For PSVM, the number of dividing points is set as 9, as Li, Artemiou and Li (2011) recommend a larger number is preferable. For a given set of dividing points, two ways to dichotomize the response are considered in  Li, Artemiou and Li (2011), ``left versus right'' (LVR) and ``one versus another''. We adopt the LVR scheme in our simulations. For PQR, we follow Wang, Shin and Wu (2018) and set the number of quantile levels to be 9, which leads to 10 slices. For PALS, we set $\uptau_k=k/10$ for $k=1,\ldots,9$. To evaluate the performance of each estimator $\hat{\B}$, we report
\begin{align}\label{delta}
\Delta=\|\bP_\B-\bP_{\hat{\B}}\|_F,
\end{align}
where $\bP_\A$ denotes the orthogonal projection onto $\mathrm{span}(\A)$, and $\|\cdot\|_F$
 is the matrix Frobenius norm. Smaller $\Delta$ value means more accurate estimation.

For the choice of the tuning parameter $\lambda$, PSVM, PQR and PALS seem to be not overly sensitive. We try $\lambda=0.1,1,10,100$, and report the best results that a fixed $\lambda$ can achieve. In addition, we propose a variable $\lambda$ scheme for PALS so that one can use different $\lambda$ values across repetitions. Specifically, denote $\hat{\B}_\lambda$ as the PALS estimator for a specific $\lambda$.  We choose $\lambda$ such that the squared sample distance correlation (Sz{\'e}kely, Rizzo and Bakirov, 2007) between 
$Y$ and $\hat{\B}^\top_\lambda \X$ is maximized. We refer to this method as DC-PALS. 

\begin{table}[t]\label{t1}
	\scalebox{.9}{\begin{tabular}{ccccccccc}
		\hline
model & $p$ & SIR & SAVE & DR & PSVM & PQR & PALS & DC-PALS \\
	\hline
\multirow{6}{*}{I} & \multirow{2}{*}{10} &1.552 & 1.779 & 1.702 &  1.482 & 1.530 & 1.424 & 1.454 \\
& & (0.016) & (0.013) & (0.017) & (0.019) & (0.018)& (0.020)& (0.018) \\
\cline{2-9}
& \multirow{2}{*}{15} & 1.698 & 1.854 & 1.800 & 1.643 & 1.643 & 1.599 & 1.606 \\
& &(0.013) & (0.010) & (0.011) & (0.011) & (0.012) & (0.014) & (0.014) \\
\cline{2-9}
& \multirow{2}{*}{20} & 1.770 & 1.914 & 1.880 & 1.712 & 1.722 & 1.672 & 1.681 \\
& & (0.010) & (0.006) & (0.008) & (0.011) & (0.010) & (0.012) & (0.011) \\
\hline

\multirow{6}{*}{II} & \multirow{2}{*}{10} &1.427 & 1.617 & 1.445 & 1.439 & 1.424 & 1.410 & 1.424 \\
& & (0.006) & (0.014) & (0.007) & (0.005) & (0.007) & (0.008) & (0.006) \\
\cline{2-9}
& \multirow{2}{*}{15} & 1.494 & 1.916 & 1.526 & 1.501 & 1.480 & 1.467 & 1.470 \\
& & (0.005) & (0.007) & (0.007) & (0.005) & (0.005) & (0.005) & (0.005) \\
\cline{2-9}
& \multirow{2}{*}{20} & 1.521 & 1.947 & 1.584 & 1.538 & 1.511 & 1.505 & 1.503 \\
& & (0.005) & (0.004) & (0.007) & (0.006) &(0.005) & (0.005) & (0.006) \\
\hline

\multirow{6}{*}{III} & \multirow{2}{*}{10} & 1.331 & 1.487 & 1.383 & 1.335 & 1.306 & 1.266 & 1.299 \\
& & (0.013) &( 0.011) & (0.009) & (0.013) & (0.016) & (0.016) & (0.017) \\
\cline{2-9}
& \multirow{2}{*}{15} & 1.413 & 1.843 & 1.429 & 1.420 & 1.360 & 1.331 & 1.378 \\
& & (0.008) & (0.012) &(0.007) & (0.008) & (0.011) & (0.014) & (0.009) \\
\cline{2-9}
& \multirow{2}{*}{20} & 1.452 & 1.924 & 1.485 & 1.458 & 1.416 & 1.408 & 1.419 \\
& & (0.006) & (0.006) & (0.006) & (0.007) & (0.008) & (0.008) & (0.007) \\
\hline
	\end{tabular}}
\caption{Results for linear sufficient dimension reduction with different $p$. The average of $\Delta$ in (\ref{delta}) and its standard error (in parenthesis) are reported based on 100 repetitions. 
}
\end{table}

First, we set $\X\sim N(\bf{0},\BSigma)$, where the element in the $i$th row and $j$th column of  $\BSigma$ is $\sigma_{i,j}=0.5^{|i-j|}$ for $i,j=1,\ldots,p$. 
We fix $n=100$, and consider  $p=10, 15, 20$. The results based on $100$ repetitions are summarized in  Table 1. We report the average of $\Delta$ in (\ref{delta}) and include its standard error in the parenthesis. 
We see that PALS with fixed $\lambda$ leads to the best result across all three models. PSVM is not as good as classical method such as SIR in model III, where heteroscedasticity is present. PQR is very competitive in models II and III, but is not as good as PSVM and PALS in model I. Furthermore, DC-PALS with variable $\lambda$ has the second best overall performance, and it is only slightly worse than PALS with fixed $\lambda$. As $p$ increases, all methods deteriorate, while PALS and DC-PALS maintain their advantage over the other methods.

\begin{table}[t]\label{t2}
	\scalebox{.9}{\begin{tabular}{ccccccccc}
		\hline
		model & case & SIR & SAVE & DR & PSVM & PQR & PALS & DC-PALS \\
		\hline
	\multirow{6}{*}{I} & \multirow{2}{*}{(i)} & 1.480 & 1.751 & 1.642 & 1.370 & 1.433 & 1.283 & 1.316 \\
	& & (0.019) & (0.014) & (0.017) & (0.023) & (0.018) & (0.021) & (0.022) \\
	\cline{2-9}
	& \multirow{2}{*}{(ii)} & 1.552 & 1.779 & 1.702 & 1.482 & 1.530 & 1.424 & 1.454 \\
	& & (0.016) & (0.013) & (0.017) & (0.019) & (0.018) & (0.020) & (0.018) \\
	\cline{2-9}
	& \multirow{2}{*}{(iii)} & 1.686 & 1.759 & 1.739 & 1.627 & 1.623 & 1.613 & 1.620 \\
	& & (0.015) & (0.013) & (0.014) & (0.015) & (0.013) & (0.014) & (0.013)\\
	\hline
		
	\multirow{6}{*}{II} & \multirow{2}{*}{(i)} &  1.383 & 1.565 & 1.364 & 1.361 & 1.358 & 1.350 & 1.357 \\
	& & (0.009) & (0.019) & (0.012) & (0.011) & (0.012) & (0.012) & (0.012) \\
	\cline{2-9}
	& \multirow{2}{*}{(ii)} & 1.427 & 1.617 & 1.445 & 1.439 & 1.424 & 1.410 & 1.424 \\
	& & (0.006) & (0.014) & (0.007) & (0.005) & (0.007) & (0.008) & (0.006) \\
	\cline{2-9}
	& \multirow{2}{*}{(iii)} & 1.404 & 1.747 & 1.447 & 1.402 & 1.400 & 1.394 & 1.393 \\
	& & (0.012) & (0.016) & (0.014) & (0.013) & (0.010) & (0.012) & (0.012) \\
	\hline	
	
	\multirow{6}{*}{III} & \multirow{2}{*}{(i)} &  1.355 & 1.451 & 1.360 & 1.330 & 1.261 & 1.203 & 1.264 \\
	& & (0.010) & (0.014) & (0.011) & (0.012) & (0.016) & (0.017) & (0.016) \\
	\cline{2-9}
	& \multirow{2}{*}{(ii)} & 1.331 & 1.487 & 1.383 & 1.335 & 1.306 & 1.266 & 1.299 \\
	& & (0.013) & (0.011) & (0.009) & (0.013) & (0.016) & (0.016) & (0.017) \\
	\cline{2-9}
	& \multirow{2}{*}{(iii)} & 1.384 & 1.581 & 1.391 & 1.365 & 1.358 & 1.355 & 1.344 \\
	& & (0.009) & (0.016) & (0.008) & (0.012) & (0.010) & (0.012) & (0.012) \\
	\hline
	\end{tabular}}
\caption{Results for linear sufficient dimension reduction with different predictor distribution. The average of $\Delta$ in (\ref{delta}) and its standard error (in parenthesis) are reported based on 100 repetitions. 
}
\end{table}

Next, we fix $n=100$, $p=10$, and consider three cases for the distribution of $\X$:  case (i),  $\X\sim N(\bf{0},I_p)$; case (ii),
$\X\sim N(\bf{0},\BSigma)$ with $\sigma_{i,j}=0.5^{|i-j|}$; and case (iii),  $X_j\sim\textrm{Uniform}(-1,1)$, $j=1,\ldots,p$, where the components of $\X$ are independent. The linear conditional mean assumption holds for cases (i) and (ii), and is no longer satisfied for case (iii).
The results based on $100$ repetitions are summarized in  Table 2.  
Compared to PALS, PQR does not work as well for model I, and PSVM is significantly worse for model III when $\X$ is normal.  
For cases (i) and (ii), all the estimators become worse when the correlation between the normal predictors increase. 
For cases (i) and (iii) with uncorrelated predictors, we see that all the methods become worse when the linear conditional mean assumption is violated. PALS and DC-PALS again have the best overall performances.

Last but not least, we list the computation time of 100 repetitions in Table 3 for PSVM, PQR and PALS when we fix $\lambda=1$ and $n=100$. We only report 
the results for model III. The other two models lead to similar results and are omitted. We see that the computation time generally increases when $p$ increases, although the increase does not seem to be significant. The predictor distribution does not seem to affect the computation time. PSVM costs the least computation time among all three methods. Although not as fast as PSVM, PALS is almost four times faster than PQR across all settings.

\begin{table}[t]\label{t3}
\centering{
  \begin{tabular}{cccccc}
		\hline
		model & case & $p$  & PSVM & PQR   & PALS \\
		\hline
		\multirow{9}{*}{III} & \multirow{3}{*}{(i)} & 10 & 3.98 &	20.91 &5.46 \\
		& & 15 & 4.34 & 21.58 & 5.66\\
		&& 20 & 4.81 & 20.79 & 5.91  \\
		\cline{2-6}
				 & \multirow{3}{*}{(ii)} & 10 & 3.86 & 20.88 & 5.57 \\
		&& 15 & 4.36 & 21.25 & 5.62 \\
		&& 20 & 4.85 & 21.57 & 5.88 \\
		\cline{2-6}
		& \multirow{3}{*}{(iii)} & 10 & 3.80 & 20.89 & 5.40 \\
		&& 15 & 4.34 & 21.10	& 5.57 \\
		&& 20 & 4.84 & 22.09 & 5.70 \\
		\hline
  \end{tabular}
}
\caption{Computation time in seconds for 100 repetitions with $\lambda = 1$.}
\end{table}

\subsection{Nonlinear sufficient dimension reduction}
For nonlinear sufficient dimension reduction, we consider the following models:

\begin{align*}
\textrm{IV: } & Y = \sqrt{\varphi_1(\X)} \log\left\{\sqrt{\varphi_1(\X)}\right\} + 0.5\varepsilon;\\
\textrm{V: } & Y= \varphi_1^2(\X) + 0.5\varphi_2(\X)\varepsilon,
\end{align*}
where $\X\sim N(\bf{0},\Ib_p)$, $\varphi_1(\X) = \sqrt{X_1^2 + X_2^2}$, $\varphi_2(\X) = \sin(X_2)$, $\varepsilon\sim N(0,0.2)$, and $\varepsilon$ is independent of $\X$.
Denote $\bvarphi (\X)$ as the basis for nonlinear sufficient dimension reduction such that $Y \indep \X | \bvarphi (\X)$. Then 
$\bvarphi (\X)=\varphi_1(\X)$ for model IV, and $\bvarphi (\X)=\{\varphi_1(\X),\varphi_2(\X)\}$ for model V.

We denote our proposal in Section 4 as kernel PALS (kPALS), and we compare it with kernel SIR (kSIR) (Wu, 2008), kernel PSVM (kPSVM) 
(Li, Artemiou and Li, 2011), and kernel PQR (kPQR) (Wang, Shin and Wu, 2018). For estimator  $\hat{\bvarphi} (\X)$, we measure its performance by the squared sample distance correlation between 
$\bvarphi (\X)$ and $\hat{\bvarphi} (\X)$ as
\begin{align}\label{dcor}
\Upsilon=\textrm{dCor}^2\left\{\bvarphi (\X),\hat{\bvarphi} (\X)\right\}.
\end{align}
Larger values of $\Upsilon$ mean better estimation. 
Similar to PSVM, PQR and PALS for linear sufficient dimension reduction, their kernel counterparts require a choice of $\lambda$. See, for example, $\lambda$ for kPALS in (\ref{obj_nonlinear2}). For $\lambda=0.1,1,10,100$, we report the results based on the best $\lambda$. Parallel to DC-PALS, we also include DC-kPALS, where $\lambda$  is chosen such that the squared sample distance correlation between 
$\hat{\bvarphi}_\lambda (\X)$ and $Y$ is maximized. 
We fix $p=10$ and set $n=100,150,200$.
From Table 4, we see that kPALS has the best performance, and DC-kPALS is a close second. 
All methods improve as $n$ increases for model IV, and only kSIR  improves as $n$ increases  for model V.
Together with previous simulation studies, we conclude that distance correlation can be a useful tool to select $\lambda$ for PALS in both linear and nonlinear sufficient dimension reduction.

\begin{table}[t]
\centering{
\begin{tabular}{ccccccc}
\hline
model & n  & kSIR & kPSVM & kPQR & kPALS & DC-kPALS \\ \hline
\multirow{6}{*}{IV}   & \multirow{2}{*}{100} & 0.523 & 0.742 & 0.750 & 0.750 & 0.750  \\ 
 &  &  (0.016) & (0.004) & (0.004) & (0.004) & (0.004)\\ \cline{2-7}
& \multirow{2}{*}{150} &  0.616 & 0.749 & 0.762 & 0.763 & 0.763  \\ 
 & &  (0.014) & (0.003) & (0.003) & (0.003) & (0.003) \\ \cline{2-7}  
& \multirow{2}{*}{200} &  0.678 & 0.756 & 0.770 & 0.772 & 0.772   \\ 
 & &  (0.010) & (0.003) & (0.003) & (0.003) & (0.003) \\ 
  \hline
  \multirow{6}{*}{V}   & \multirow{2}{*}{100} & 0.496 & 0.578 & 0.583 & 0.605 & 0.599   \\ 
 &  &  (0.008) & (0.003) & (0.004) & (0.004) & (0.003) \\ \cline{2-7}
& \multirow{2}{*}{150} & 0.537 & 0.580 & 0.581 & 0.606 & 0.602   \\ 
 & &  (0.005) & (0.003) & (0.003) & (0.003) & (0.003) \\ \cline{2-7}  
& \multirow{2}{*}{200} &  0.554 & 0.579 & 0.582 & 0.605 & 0.598  \\ 
 & &  (0.004) & (0.002) & (0.003) & (0.002) & (0.002)\\ 
  \hline
\end{tabular}}
\label{t4}
\caption{Results for nonlinear sufficient dimension reduction with different $n$. The average of $\Upsilon$ in (\ref{dcor}) and its standard error (in parenthesis) are reported based on 100 repetitions. }
\end{table}

\section{Real data analysis of the Boston housing data}\label{real_anal}

We consider Boston housing data for the real data analysis. The data is originally studied in Harrison and Rubinfeld (1978).
After removing a categorical variable and excluding the cases where the census tract bounds the Charles river, we end up with 12 predictors and 471 observations. The response is the median value of owner-occupied homes in each census tract. For a complete list of the predictors, one can refer to Table 5 of
Wang, Shin and Wu  (2018).
As suggested by Wang, Shin and Wu  (2018), we set the structural dimension to be $d=1$ and denote the estimator as $\hat{\bbeta}$.
 We apply PSVM, PQR and PALS to this data set, and report the squared sample distance correlation between $Y$ and $\hat{\bbeta}^\top\X$
for different $\lambda$. From Table 5, we see that PQR and PALS perform similarly, and both are better than PSVM. Furthermore, PQR and PALS are less sensitive to the choice of $\lambda$ than PSVM. 

\begin{table}[t]
\centering{
\begin{tabular}{ccccc}
\hline
 & $\lambda=0.1$ & $\lambda=1$  &  $\lambda=10$ & $\lambda=100$   \\
   \hline
  PSVM & 0.831  & 0.812 &  0.721 &  0.711 \\
   \hline
   PQR &  0.866 &  0.864 & 0.864 & 0.864\\
    \hline
  PALS &  0.863 & 0.863 & 0.863 & 0.864 \\
 \hline
\end{tabular}}
%\vspace{-0.2in}
\caption{The squared sample distance correlation between $Y$ and $\hat{\bbeta}^\top\X$ for the Boston housing data.}
\label{t5}
\end{table}

\section{Conclusion}
We propose PALS for linear and nonlinear sufficient dimension reduction in this paper. Our proposed method is very competitive with existing methods in the literature. 
On one hand, our proposal enjoys better estimation accuracy than SIR and PSVM, especially in the presence of heteroscedasticity. 
On the other hand, our proposal is computationally more efficient compared to PQR. Unlike PSVM where the response is dichotomized, both PQR and PALS deal with
 continuous response directly. We develop new tools for the asymptotic analysis of PALS. Specifically, Lemma 3 of Li, Artemiou and Li (2011) provides a tool to compute the derivative of an expectation of a non-Lipschitz function, and Theorem 3 of  Wang, Shin and Wu  (2018)
applied this Lemma directly without considering the continuous support of the response in PQR. This limitation is addressed in Lemma 1 and Theorem 5 of the Appendix, where Lemma 3 of Li, Artemiou and Li (2011) is adapted for continuous response. 

We consider a fixed set of expectile levels in this paper. Although our experience indicates that the performance of PALS is not sensitive to the choice of expectile levels, choosing an optimal set of expectile levels is worth further investigation. Kim, Wu and Shin (2019) develop quantile-slicing for sufficient dimension reduction, and  expectile-slicing for sufficient dimension reduction may be an interesting research direction.

%\appendix
\section*{Appendix A}

\noindent {\bf Proof  of Theorem \ref{thm_linearity}}.
	We assume without loss of generality that $E(\X) = \bm{0}$. Note that $\var(\bm{\beta}^\top\X)=\bm{\beta}^\top\BSigma\bm{\beta}$. Then  (\ref{obj_pop}) becomes
\begin{align}\label{a1}
		L_\uptau(\alpha, \bm{\beta}) =  \var(\bm{\beta}^\top   \X) + \lambda E\{ \rho_\uptau (Y-\alpha-\bm{\beta}^\top\X ) \}. 
		\end{align}	
	The first term on the right hand side of 	(\ref{a1}) satisfies
\begin{align}\label{a2}
	\var(\bm{\beta}^\top   \X)\geq
	\var\{E(\bm{\beta}^\top   \X | \B^\top   \X)  \}.
		\end{align}	
	%	Let's first consider $var(\bm{\beta^\top   X})$. Using the conditional variance formula, it can be decomposed as
%	\begin{align*}
%	var(\bm{\beta^\top   X}) &= E\left[var(\bm{\beta^\top   X} | \bm{B^\top   X})\right] + var\left[E(\bm{\beta^\top   X} | \bm{B^\top   X})\right]\\
%	&\geq
%	var\left[E(\bm{\beta^\top   X} | \bm{B^\top   X}) \right]
%	\end{align*}
%	since $E\left[var(\bm{\beta^\top   X} | \bm{B^\top   X}) \right] \geq \bm{0}$.\\
The second  term on the right hand side of 	(\ref{a1}) satisfies
	\begin{align}\label{a3}
	\begin{split}
	E\{ \rho_\uptau(Y-\alpha-\bm{\beta}^\top  \X) \} &=
	E[ E\{ \rho_\uptau(Y-\alpha-\bm{\beta}^\top\X) | \B^\top \X,Y \} ] \\
	& \geq
	E[ \rho_\uptau\{ E(Y-\alpha-\bm{\beta}^\top\X) | \B^\top \X,Y \} ] \\
	&=
	E [\rho_\uptau \{ Y-\alpha-E(\bm{\beta}^\top \X|\B^\top \X ) \}],
	\end{split}
	\end{align}
	where the inequlality is due to the convexity of $\rho_\uptau$, and the last equality is due to the conditional independence (\ref{SDR}).
The assumption that 	$E(\X|\B^\top \X )$ is linear in $\B^\top \X$ implies that
\begin{align}\label{a4}
E(\X|\B^\top \X )=\BSigma (\B^\top\BSigma\B ) ^{-1} \B^\top \X.
\end{align}
(\ref{a1}), (\ref{a2}), (\ref{a3}) and (\ref{a4}) together imply that
$$L_\uptau(\alpha, \bm{\beta})\geq L_\uptau(\alpha, \tilde{\bm{\beta}}) \mbox{ with }\tilde{\bm{\beta}}=\B (\B^\top\BSigma\B ) ^{-1} \BSigma \bm{\beta}.$$
Thus the minimizer $\bm{\beta}_{0,\uptau}$ must satisfy $\bm{\beta}_{0,\uptau}=\B (\B^\top\BSigma\B ) ^{-1} \BSigma \bm{\beta}_{0,\uptau}\in \mathrm{span}(\B)=\spc$. \eop

%  Since $\varrho_\uptau(.)$ is convex, the second inequality holds by Jensen's inequality for conditional expectations. This implies that the minimum is always achieved at $E(\bm{\beta^\top {X}}|\bm{B^\top   X})$ even for a possibly non-unique estimate.  Therefore by the linearity condition $E(\bm{\beta^\top   {X}}|\bm{B^\top   X}) \in \spc$, where $E(\bm{\beta^\top   {X}}|\bm{B^\top   X}) = \bm{\beta}^\top   P_{\bm{B}}^\top  (\bm{\Sigma})\bm{X}$.\\
%	\indent Now, suppose $\exists \bm{\beta}_1 \notin \spc$, then $var(\bm{\beta}_1 ^\top   \bm{X} | \bm{B^\top    X} )$ must be strictly greater than $0$ which is a contradiction. Hence $\bm{\beta}_1$ cannot be a minimizer.
%\end{proof}

\smallskip

\noindent {\bf Proof  of Corollary \ref{co}}. The proof follows directly from Theorem \ref{thm_linearity} and is omitted. \eop

\smallskip

\noindent  {\bf Proof  of Theorem \ref{thm_duality}}. Denote $\bm{\xi_{+}}=(\xi_{1+},\ldots,\xi_{n+})^\top$, 
$\bm{\xi_{-}}=(\xi_{1-},\ldots,\xi_{n-})^\top$, $\bm{u}=(u_1,\ldots,u_n)^\top$,  $\bm{v}=(v_1,\ldots,v_n)^\top$,
$\bm{a}=(a_1,\ldots,a_n)^\top$, and $\bm{\eta}=(\eta_1,\ldots,\eta_n)^\top$. 
Denote $L^*(\alpha,\bm{\theta},\bm{\xi_{+}},\bm{\xi_{-}},\bm{u},\bm{v},\bm{a},\bm{\eta})$ as the Lagrangian of the primal optimization problem (\ref{opt_primal}) and abbreviate it as $L^*$. Then we have
%\begin{proof}[Proof of Theorem \ref{thm_duality}] \quad \\
%	Define $\xi_{i+} = (Y_i - \alpha- \bm{\theta}^\top   \bm{Z}_i)_{+}$
%	and $\xi_{i-} = ( \alpha+ \bm{\theta}^\top   \bm{Z}_i-Y_i)_{+}$, where $c_{+}=max(0,c)$. Then (\ref{obj_sam2}) leads to primal optimization problem
%	\begin{align}
%	\underset{\alpha, \bm{\theta}}{\operatorname{min}} \hspace{0.11in} \bm{\theta}^\top  \bm{\theta} + \lambda\uptau \sum_{i=1}^n \xi_{i+}^2 + \lambda(1-\uptau) \sum_{i=1}^n \xi_{i-}^2
%	\label{opt primal}
%	\end{align}
%	subject to $\xi_{i+} \geq 0$, $\xi_{i-}\geq 0$,
%	$\xi_{i+} \geq Y_i - \alpha- \bm{\theta}^\top   \bm{Z}_i$, and $\xi_{i-} \geq \alpha+ \bm{\theta}^\top   \bm{Z}_i-Y_i$. \footnote{Recall that $1/n$ is ignored. Hence, $\lambda$ is used instead of $\cfrac{\lambda}{n}$}. \\
%	
%	The Lagrangian of the primal is given by \\

	\begin{equation}
	\begin{aligned}
	L^*=  & \hspace{0.11in} \bm{\theta}^\top  \bm{\theta} + \tlambda\uptau \sum_{i=1}^n \xi_{i+}^2 + \tlambda(1-\uptau) \sum_{i=1}^n \xi_{i-}^2  - \sum_{i=1}^n u_i\xi_{i+} - \sum_{i=1}^n v_i\xi_{i-} \\
	& + \sum_{i=1}^n a_i(Y_i - \alpha - \bm{\theta}^\top   \Z_i-\xi_{i+}) +
	\sum_{i=1}^n \eta_i(-Y_i + \alpha + \bm{\theta}^\top   \Z_i-\xi_{i-}),
	\end{aligned}
	\label{Langrangian}
	\end{equation}
	where $u_i\geq 0$, $v_i\geq 0$, $a_i \geq 0$, and $\eta_i \geq 0$ for all $i$.
%	where $a, \bm{\eta, u,v} \in \mat hbb{R}^n$ and $a_i, \eta_i, u_i,v_i > 0$ . \\
%	\newline
	Take partial derivatives of (\ref{Langrangian}) and set them to be zero. We get
	\begin{align}\label{b1}
	\begin{cases}
	&{\partial L^*}/{\partial \bm{\theta}}=2\bm{\theta} -\sum_{i=1}^n (a_i - \eta_i)\Z_i = \bf{0}\\
&	{\partial L^*}/{\partial \alpha}=\sum_{i=1}^n (\eta_i - a_i)=0\\
& {\partial L^*}/{\partial \xi_{i+}}=2\tlambda\uptau \xi_{i+} - u_i - a_i  = 0\\
& {\partial L^*}/{\partial \xi_{i-}}= 2\tlambda(1-\uptau)\xi_{i-} -v_i-\eta_i=0
\end{cases}
\end{align}	
 Assume $u_i>0$ for a particular $i$.  The  Karush Kuhn Tucker (KKT) conditions state that
$u_i \xi_{i+} = 0$ for all $i$. Then we must have $\xi_{i+}=0$ from KKT. On the other hand, we have $\xi_{i+}=(u_i + a_i)/(2\tlambda\uptau)$ from the third equation of (\ref{b1}), which leads to $\xi_{i+}>0$ because  $u_i>0$, $a_i\ge 0$, $\tlambda>0$ and $\uptau>0$. This contradiction guarantees that $u_i=0$ for all $i$. Thus we have
 \begin{align}\label{b2}
 \xi_{i+}= \frac{a_i}{2\tlambda\uptau} \mbox{ for all }i.
 \end{align}	
 Similarly, from the fourth equation of (\ref{b1}) and the KKT condition, we have $v_i=0$ for all $i$ and
 \begin{align}\label{b3}
 \xi_{i-}= \frac{\eta_i}{2\tlambda(1-\uptau)} \mbox{ for all }i.
 \end{align}		
 Furthermore, the first equation of (\ref{b1}) leads to
  \begin{align}\label{b4}
  \bm{\theta} = \frac{1}{2}\sum_{i=1}^n (a_i - \eta_i)\Z_i 
   \end{align}	
	By complementary slackness, we have
		\begin{align}\label{b5}
		\begin{split}
	& u_i \xi_{i+} = 0,  
	 a_i(Y_i - \alpha - \bm{\theta}^\top  \Z_i-\xi_{i+}) = 0,\\
	 & v_i \xi_{i-} = 0,
	\eta_i(-Y_i + \alpha + \bm{\theta}^\top  \Z_i-\xi_{i-}) = 0 \mbox{ for all }i.
	\end{split}
	\end{align}
Plug (\ref{b2}), (\ref{b3}), (\ref{b4}) and (\ref{b5}) into (\ref{Langrangian}), and we get the objective function in the dual optimization problem (\ref{opt_dual}). The constraints for the dual problem are $a_i \geq 0$ and $\eta_i \geq 0$ for all $i$. The second equation of (\ref{b1}) leads to
the constraint that $a_i=\eta_i$ for all $i$. Equation  (\ref{b4}) leads to (\ref{pd}), which connects the solution of the dual problem to the primal problem. \eop

\section*{Appendix B}

We provide the proof of  Theorem \ref{normal} in this section.
The following notations are needed.
% before we state the assumptions needed for
Without loss of generality, assume $E(\X)=\bf{0}$.
Denote $\BXi=(\X^\top,Y)^\top$, $\tX=(1,\X^\top)^\top$ and $\tbeta=(\alpha,\bm{\beta}^\top)^\top$. Let $\tSigma\in\R^{(p+1)\times(p+1)}$ be a block diagonal matrix such that the block diagonal elements of $\tSigma$ are $0$ and $\BSigma$.
Then $L_\uptau(\alpha,\bm{\beta})$
% =  \bm{\beta}^\top   \bm{\bm{\Sigma}} \bm{\beta} + \lambda E\left[ \rho_\uptau\left\{Y-\alpha-\bm{\beta}^\top  (\X-\bm{\mu} )\right\} \right].
in (\ref{obj_pop}) becomes $E\{\ell_\uptau(\tbeta,\BXi)\}$, where
\begin{align}\label{ell}
\ell_\uptau(\tbeta,\BXi)=\tbeta^\top\tSigma\tbeta+\lambda  \rho_\uptau (Y-\tbeta^\top \tX). 
\end{align}
Let $D_{\tbeta}$ be the $(p+1)$-dimensional column vector of differential operators $(\partial/\partial \alpha,\partial/\partial \beta_1,\ldots,\partial/\partial\beta_p)^\top$. The next result gives the gradient of $E\{\ell_\uptau(\tbeta,\BXi)\}$.

\begin{theorem}\label{d}
Suppose for any $y$, the distribution of $\X|Y=y$ is dominated by the Lebesgue measure, $E(Y^2)<\infty$  and $E(\|\X\|^2)<\infty$.
Then
\begin{align}\label{dd}
%\begin{split}
%D_{\tbeta}E\{\ell_\uptau(\tbeta,\BXi)\}&=(0,2\bbeta^\top\BSigma)^\top-2\lambda\uptau E(\xi_{+}\tX)\\
%&\hspace{.2in} +2\lambda E\{\tX(\uptau\xi+\xi_{-})\mI(\xi<0)\}
%\end{split}
D_{\tbeta}E\{\ell_\uptau(\tbeta,\BXi)\}=(0,2\bbeta^\top\BSigma)^\top-2\lambda E\{\uptau \xi_{+}\tX-(\uptau\xi+\xi_{-})\mI(\xi<0)\tX\},
\end{align}
where $\xi=Y-\tbeta^\top \tX$, $\xi_+=max(\xi,0)$, $\xi_-=max(-\xi,0)$, and $\mI(\cdot)$ is the indicator function.
\end{theorem}

\noindent  {\bf Proof}. Denote $\ell_\uptau^*(\tbeta,\BXi)=\rho_\uptau (Y-\tbeta^\top \tX)$. It is easy to check that 
\begin{align}\label{d1}
\ell_\uptau^*(\tbeta,\BXi)=\uptau \xi_+^2+(1-\uptau)\xi_-^2.
\end{align}
Note that $D_{\tbeta}\xi_+=-\{1-\mI(\xi<0)\}\tX$. It follows that
\begin{align}\label{d2}
D_{\tbeta}\xi_+^2=-2\xi_+\{1-\mI(\xi<0)\}\tX.
\end{align}
Similarly from $D_{\tbeta}\xi_-=\mI(\xi<0)\tX$, we have
\begin{align}\label{d3}
D_{\tbeta}\xi_-^2=2\xi_-\mI(\xi<0)\tX.
\end{align}
Plug (\ref{d2}) and (\ref{d3}) into (\ref{d1}). After taking derivatives, we get
\begin{align}\label{d4} 
\begin{split}
D_{\tbeta} \ell_\uptau^*(\tbeta,\BXi)&=-2\uptau \xi_+\{1-\mI(\xi<0)\}\tX+2(1-\uptau)\xi_-\mI(\xi<0)\tX\\
&= -2\uptau \xi_+\tX+2 \uptau (\xi_+-\xi_-)\mI(\xi<0)\tX+2\xi_-\mI(\xi<0)\tX \\
&=-2\uptau \xi_+\tX+2(\uptau\xi+\xi_{-})\mI(\xi<0)\tX.
\end{split}
\end{align}
(\ref{d4}) and 
(\ref{ell}) together imply that 
\begin{align}\label{d7}
E \{D_{\tbeta} \ell_\uptau(\tbeta,\BXi)\}=(0,2\bbeta^\top\BSigma)^\top-2\lambda E\{\uptau \xi_{+}\tX-(\uptau\xi+\xi_{-})\mI(\xi<0)\tX\}.
\end{align}
%It remain to show $D_{\tbeta}E\{\ell_\uptau(\tbeta,\BXi)\}=E \{D_{\tbeta} \ell_\uptau(\tbeta,\BXi)\}$. From Lemma 2 of Li, Artemiou and Li (2011), we need  to satisfy the Lipschitz condition with respect to $\tbeta$.  It suffices to check the Lipschitz condition for 
%$\ell_\uptau^*(\tbeta,\BXi)$. 

Let $\Theta$ be the support of $\tbeta$.  For $\tbeta_1=(\alpha_1,\bm{\beta}_1^\top)^\top\in\Theta$  and  $\tbeta_2=(\alpha_2,\bm{\beta}_2^\top)^\top\in\Theta$,   we have
\begin{align}\label{d5} \begin{split}
\ell_\uptau^*(\tbeta_1,\BXi)-\ell_\uptau^*(\tbeta_2,\BXi)&=\uptau\{(Y-\tbeta_1^\top \tX)_+^2-(Y-\tbeta_2^\top \tX)_+^2\}\\
 &\hspace{.1in}+ (1-\uptau)\{(\tbeta_1^\top \tX-Y)_+^2-(\tbeta_2^\top \tX-Y)_+^2\}.
 \end{split}
\end{align}
Note that $u_+-v_+\leq |u-v|$ and $u_+ + v_+\leq|u|+|v| $. Then 
\begin{align*}
%\label{d6} 
&(Y-\tbeta_1^\top \tX)_+^2-(Y-\tbeta_2^\top \tX)_+^2
\leq  |(\tbeta_1-\tbeta_2)^\top \tX| (|Y-\tbeta_1^\top \tX|+|Y-\tbeta_2^\top \tX|)\\
 &\hspace{.1in} \leq (1+\|\X^2\|)^{1/2}(|Y-\tbeta_1^\top \tX|+|Y-\tbeta_2^\top \tX|)\|\tbeta_1-\tbeta_2\|<c\|\tbeta_1-\tbeta_2\|
\end{align*}
for some constant $c<\infty$. The last inequality is due to the assumption that  $E(Y^2)<\infty$  and $E(\|\X\|^2)<\infty$. Thus the first term on the right hand side of (\ref{d5}) satisfies the Lipschitz condition with respect to $\tbeta$. Similarly, one can show the second term on the right hand side of (\ref{d5}) also satisfies the Lipschitz condition.  Together, we know  $\ell_\uptau(\tbeta,\BXi)$ satisfies the Lipschitz condition with respect to $\tbeta$.
From Lemma 2 of Li, Artemiou and Li (2011), we have
\begin{align}\label{d6}
D_{\tbeta}E\{\ell_\uptau(\tbeta,\BXi)\}=E \{D_{\tbeta} \ell_\uptau(\tbeta,\BXi)\}.
\end{align}
(\ref{d7}) and (\ref{d6}) together lead to the desired result.   
  \eop

\smallskip

The next lemma is used to compute the derivative of an expectation of a non-Lipschitz function. Let $D_{\epsilon=0}$ denote the operation of first taking derivative with respect to $\epsilon$ and then evaluating the derivative at $\epsilon=0$.

  \begin{lemma}\label{uvw}
  Suppose $U$, $V$ and $W$ are random variables, and $\h(u,v,w)$ is a measurable $\R^k$-valued function. Suppose, moreover,
  \begin{itemize}
  \item[1.] the joint distribution of $(U,V,W)$ is dominated by the Lebesgue measure;
  \item[2.] for any $(v,w)$, the function $u\mapsto \h(u,v,w)f_{U|V,W}(u|v,w)$ is continuous, where  $f_{U|V,W}$ denotes the conditional density of $U$ given $(V,W)$;
  \item[3.] for each component $h_i(u,v,w)$ of $\h(u,v,w)$, there is a function $c_i(v,w)\ge 0$ such that
  $$h_i(u,v,w)f_{U|V,W}(u|v,w)\le c_i(v,w) \mbox{ and } E\{c_i(V,W)\}<\infty.$$
  \end{itemize}
  Then, for any constant $a$, the function $$\epsilon\mapsto E\{\h(U,V,W)\mI(W+U+\epsilon V<a+\epsilon \eta)\}$$
  is differentiable at $\epsilon=0$ with derivative 
  \begin{align}\label{e}
  \begin{split}
&D_{\epsilon=0} E\{\h(U,V,W)\mI(W+U+\epsilon V<a+\epsilon \eta)\}\\
 &  =E_W[ f_{U|W}(a-W|W) E_V\{(\eta-V)h_i(a-W,V,W)|U=a-W,W\}],
\end{split}
\end{align}
where $E(\cdot)$ is with respect to the joint distribution of $(U,V,W)$, $E_W(\cdot)$ is with respect to the marginal distribution of $W$, and
$E_V(\cdot|U=a-W,W)$ is with respect to the conditional distribution $f_{V|U,W}(v|a-w,w)$.
  \end{lemma}
  
\noindent  {\bf Proof}. By the mean value theorem and assumptions 2 and 3, there is a $\delta\in(0,\epsilon)$  such that 
\begin{align*}
&  \left| \epsilon^{-1}\int_{a-w}^{a-w+\epsilon(\eta-v)} h_i(u,v,w)f_{U|V,W}(u|v,w)du\right|\\
 &\hspace{.2in}  =|h_i(a-w+\delta(\eta-v),v,w)f_{U|V,W}(a-w+\delta(\eta-v)|v,w)|\leq c_i(v,w)
  \end{align*}
  By the dominated convergence theorem, we have
  \begin{align*}
&  \lim_{\epsilon\rightarrow 0}\int\int\left\{  \epsilon^{-1}\int_{a-w}^{a-w+\epsilon(\eta-v)} h_i(u,v,w)f_{U|V,W}(u|v,w)du\right\}f_{V,W}(v,w)dvdw\\
&=\int\int \lim_{\epsilon\rightarrow 0}\left\{  \epsilon^{-1}\int_{a-w}^{a-w+\epsilon(\eta-v)} h_i(u,v,w)f_{U|V,W}(u|v,w)du\right\}f_{V,W}(v,w)dvdw\\
&= \int\int (\eta-v) h_i(a-w,v,w)f_{U|V,W}(a-w|v,w)f_{V,W}(v,w)dvdw\\
&= \int \left\{ f_W (w)f_{U|W}(a-w|w)\int (\eta-v) h_i(a-w,v,w)f_{V|U,W}(v|a-w,w)dv\right\}dw\\
&= E_W[ f_{U|W}(a-W|W) E_V\{(\eta-V)h_i(a-W,V,W)|U=a-W,W\}].
    \end{align*}
    Here $f_{V|U,W}$ and $f_{U|W}$ are conditional density functions, and $f_W$ denotes the marginal density of $W$.    
    \eop

\smallskip

We now present the hessian matrix of the PALS objective function (\ref{obj_pop}) in the next Theorem. 
%Lemma \ref{uvw} allows us to prov Although 

\begin{theorem}\label{h}
Suppose $\X$ has a convex and open support and its conditional distribution given $Y=y$ for any $y\in\R$ is dominated by the Lebesgue measure. Suppose, moreover,
  \begin{itemize}
  \item[1.] for any linearly independent $\bbeta,\bdelta\in\R^p$ and any $y\in\R$, the following function is continuous
  $$u\mapsto E\{\tX (\uptau\xi+\xi_{-})|-\bbeta^\top \X=u,\bdelta^\top \X=v,Y=y\}f_{-\bbeta^\top \X|\bdelta^\top \X,Y}(u|v,y);$$
  \item[2.] for any $i=1,\ldots,p$ and  any $y\in\R$, there is a nonnegative function $c_i(v,y)$ with $E\{c_i(V,Y)\}<\infty$ such that
  $$E\{X_i (\uptau\xi+\xi_{-})|-\bbeta^\top \X=u,\bdelta^\top \X=v,Y=y\}f_{-\bbeta^\top \X|\bdelta^\top \X,Y}(u|v,y)\le c_i(v,y);$$
  \item[3.] there is a nonnegative function $c_0(v,y)$  with  $E\{c_0(V,Y)\}<\infty$ such that
    $$f_{-\bbeta^\top \X|\bdelta^\top \X,Y}(u|v,y)\le c_0(v,y).$$
    \end{itemize}
    Then the function $\tbeta\mapsto D_{\tbeta}E\{\ell_\uptau(\tbeta,\BXi)\}$ is differential in all directions with derivative matrix 
      \begin{align*}
\bH_\uptau=& 2 \text{diag}(0,\BSigma)+2\lambda \uptau E[\{1-\mI(\xi<0)\}\tX \tX^\top]\\
 &\hspace{.2in} +2\lambda E_Y[ f_{-\bbeta^\top \X|Y}(\alpha-Y|Y) E_{\X}\{(\uptau\xi+\xi_{-})\tX \tX^\top|-\bbeta^\top \X=\alpha-Y,Y\}],
    \end{align*}
    where $E(\cdot)$ is  with respect to the joint distribution of $(\X,Y)$, $E_Y(\cdot)$ is  with respect to the marginal distribution of $Y$, and
    $E_{\X}(\cdot|-\bbeta^\top \X=\alpha-Y,Y)$ is with respect to the conditional distribution  $f_{\X|-\bbeta^\top \X,Y}(\bm{x}|\alpha-y,y)$. Furthermore, if the function $\tbeta\mapsto \uptau E[\{1-\mI(\xi<0)\}\tX \tX^\top]
  + E_Y[ f_{-\bbeta^\top \X|Y}(\alpha-Y|Y) E_{\X}\{(\uptau\xi+\xi_{-})\tX \tX^\top|-\bbeta^\top \X=\alpha-Y,Y\}] $ is continuous, then 
  $D_{\tbeta}E\{\ell_\uptau(\tbeta,\BXi)\}$ is jointly differentiable with respect to $\tbeta$.
\end{theorem}

\noindent  {\bf Proof}. Recall that  $\xi=Y-\bbeta^\top \X-\alpha$ and $\xi_-=max(-\xi,0)$.  First, we verify the directional differentiability of the function  $\tbeta\mapsto E\{(\uptau\xi+\xi_{-})\mI(\xi<0)\tX\}$.
For $\bdelta\in\R^p$ and $\eta\in\R$, the directional derivative along $(\eta,\bdelta^\top)^\top$ is the derivative of the following function with respect to $\epsilon$ at $\epsilon=0$,
\begin{align*}
&E\{\tX (\uptau\xi+\xi_{-}) \mI(Y-\bbeta^\top \X-\alpha+\epsilon \bdelta^\top \X<\epsilon \eta)\}\\
&\hspace{.2in} =E[ E \{\tX (\uptau\xi+\xi_{-}) |Y, \bbeta^\top \X,\bdelta^\top \X \} \mI(Y-\bbeta^\top \X+\epsilon \bdelta^\top \X<\alpha+\epsilon \eta)]
\end{align*}
Let $W=Y$, $U=-\bbeta^\top \X$, $V=\bdelta^\top \X$, $\h(U,V,W)=E\{\tX (\uptau\xi+\xi_{-})|U,V,W\}$, and $a=\alpha$. By (\ref{e}) in Lemma \ref{uvw}, the derivative above is
\begin{align*}
E_Y[ f_{-\bbeta^\top \X|Y}(\alpha-Y|Y) E_{\bdelta^\top \X}\{\tX (\uptau\xi+\xi_{-})(\eta-\bdelta^\top \X)|-\bbeta^\top \X=\alpha-Y,Y\}].
\end{align*}
Since this holds for all $(\eta,\bdelta^\top)^\top$, the function  $\tbeta\mapsto E\{(\uptau\xi+\xi_{-})\mI(\xi<0)\tX\}$ is directionally differentiable with derivative matrix
\begin{align}\label{h1}
E_Y[ f_{-\bbeta^\top \X|Y}(\alpha-Y|Y) E_{\X}\{(\uptau\xi+\xi_{-})\tX \tX^\top|-\bbeta^\top \X=\alpha-Y,Y\}].
\end{align}
On the other hand, it is easy to see 
\begin{align}\label{h2}
D_{\tbeta} (0,\bbeta^\top\BSigma)^\top=\text{diag}(0,\BSigma)
\end{align}
 and 
\begin{align}\label{h3}
D_{\tbeta}E(\xi_{+}\tX)=-E[\{1-\mI(\xi<0)\}\tX \tX^\top].
\end{align}
Plug (\ref{h1}), (\ref{h2}) and (\ref{h3})  into (\ref{dd}), and we get the desired result. \eop

\smallskip

Let $\tbeta_{0,\uptau}=(\alpha_{0,\uptau},\bm{\beta}_{0,\uptau}^\top)^\top$ from minimizing $L_\uptau(\alpha,\bm{\beta})$
in (\ref{obj_pop}). Let $\hat{\tbeta}_{0,\uptau}=(\hat{\alpha}_{0,\uptau},\hat{\bm{\beta}}_{0,\uptau}^\top)^\top$ from minimizing $\hat{L}_\uptau(\alpha,\bm{\beta})$ in
(\ref{obj_sam1}). The next result gives the influence function of PALS. Its proof follows Theorem 5.23 of Van der Vaart (2000), and is thus omitted. 

\begin{theorem}\label{inf}
If the conditions in Theorem \ref{d} and Theorem \ref{h} are satisfied, then 
\begin{align*}
\hat{\tbeta}_{0,\uptau}=& \tbeta_{0,\uptau}-\bH_{0,\uptau}^{-1}\{
(0,2\bbeta_{0,\uptau}^\top\BSigma)^\top-2\lambda\uptau E_n(\xi_{0,\uptau}^+\tX)\\
&\hspace{.02in} +2\lambda E_n\{\tX(\uptau\xi_{0,\uptau}+\xi_{0,\uptau}^-)\mI(\xi_{0,\uptau}<0)\}+o_P(n^{-1/2}),
\end{align*}
where  $\xi_{0,\uptau}=Y-\tbeta_{0,\uptau}^\top \tX$, $\xi_{0,\uptau}^+=max(\xi_{0,\uptau},0)$, $\xi_{0,\uptau}^-=max(-\xi_{0,\uptau},0)$, $E_n(\cdot)$ is with respect to the empirical distribution of $(\X,Y)$, and
      \begin{align*}
&\bH_{0,\uptau}= 2 \text{diag}(0,\BSigma)+2\lambda \uptau E[\{1-\mI(\xi_{0,\uptau}<0)\}\tX \tX^\top]\\
 &\hspace{.2in} +2\lambda E_Y[ f_{-\bbeta_{0,\uptau}^\top \X|Y}(\alpha-Y|Y) E_{\X}\{(\uptau\xi_{0,\uptau}+\xi_{0,\uptau}^{-})\tX \tX^\top|-\bbeta_{0,\uptau}^\top \X=\alpha-Y,Y\}].
    \end{align*}
\end{theorem}

\smallskip

\noindent {\bf Proof of  Theorem \ref{normal}}.
Consider $0<\uptau_1<\ldots<\uptau_K<1$. For any $k=1,\ldots,K$, let $\F_k$ be the last $p$ rows of $\bH_{0,\uptau_k}^{-1}$. Denote
\begin{align*}
%\s_k(\tbeta_{0,\uptau_k},\BXi)=&\F_k\{
%(0,2\bbeta_{0,\uptau_k}^\top\BSigma)^\top-2\lambda\uptau_k \xi_{0,\uptau_k}^{+}\tX\\
%&\hspace{.02in} +2\lambda\tX(\uptau_k\xi_{0,\uptau_k}+\xi_{0,\uptau_k}^{-})\mI(\xi_{0,\uptau_k}<0),
\s_k(\tbeta_{0,\uptau_k},\BXi)=\F_k\{
(0,2\bbeta_{0,\uptau_k}^\top\BSigma)^\top-2\lambda\uptau_k \xi_{0,\uptau_k}^{+}\tX
 +2\lambda\tX(\uptau_k\xi_{0,\uptau_k}+\xi_{0,\uptau_k}^{-})\mI(\xi_{0,\uptau_k}<0),
\end{align*}
%\begin{align}\label{f1}
%\begin{split}
%\s_k(\tbeta_{0,\uptau_k},\BXi_i)=&\F_k\{
%(0,2\bbeta_{0,\uptau_k}^\top\BSigma)^\top-2\lambda\uptau_k \xi_{0,\uptau_k}^{i,+}\tX_i\\
%&\hspace{.02in} +2\lambda\tX_i(\uptau_k\xi_{0,\uptau_k}^i+\xi_{0,\uptau_k}^{i,-})\mI(\xi_{0,\uptau_k}^i<0),
%\end{split}
%\end{align}
%where $\xi_{0,\uptau_k}^i=Y_i-\tbeta_{0,\uptau_k}^\top \tX_i$, $\xi_{0,\uptau_k}^{i,+}=max(\xi_{0,\uptau_k}^i,0)$, and $\xi_{0,\uptau_k}^{i,-}=max(-\xi_{0,\uptau_k}^i,0)$.  
For a matrix $\A\in\R^{r_1\times r_2}$, let $\K_{r_1,r_2}\in\R^{r_1 r_2\times r_1r_2}$ be the commutation matrix defined by the relation 
$\K_{r_1,r_2}\vect(\A)=\vect(\A^\top)$. Define $\BOmega_1=\Ib_{p^2}+\K_{p,p}$,
$$\BOmega_2=\sum_{k=1}^K  \sum_{t=1}^K (\bm{\beta}_{0,\uptau_k} \bm{\beta}_{0,\uptau_t}^\top) \otimes  E\{\s_k(\tbeta_{0,\uptau_k},\BXi)
\s_t^\top(\tbeta_{0,\uptau_t},\BXi)\} ,$$
and  $\BOmega=\BOmega_1 \BOmega_2 \BOmega_1$. The result then follows from Theorem \ref{inf} above and  Theorem 7 in  Li, Artemiou and Li (2011). \eop


\begin{thebibliography}{99}

\bibitem{Abdous1995}
Abdous, B. and  Remillard, B. (1995) Relating quantiles and expectiles under weighted-symmetry. {\it Annals of the Institute of Statistical Mathematics},  {\bf 47}, 371--384.

\bibitem{Artemiou2016} Artemiou, A. and Dong, Y. (2016) Sufficient dimension reduction via principal $\ell$-q support vector machine. {\it Electronic Journal of Statistics}, {\bf 10}, 783--805.

%\bibitem{}
%Bura, E. and Pfeiffer, R. (2008) On the distribution of the left singular vectors of a random matrix and its applications.{\it Statist. Probab. Lett.} {\bf 78} 2275–2280.
%\bibitem{}
%Bondell, H. D., Reich, B. J. and Wang, H. (2010). Non-crossing quantile regression curve estimation. {\it Biometrika}, {\bf 97}, 825--838.

%\bibitem{}
%Chen, J.,  Stern, M.,  Wainwright, M. and  Jordan, M. I. (2018)
%Kernel feature selection via conditional covariance minimization.  In Bengio, S.,  Fergus,  R.,  Vishwanathan, S. and  Wallach, H. (Eds), {\it Advances in Neural Information Processing Systems}, {\bf 31}.

%\bibitem{}
%Cook, R.D. (1977) Detection of influential observation in linear regression. {\it Technometrics}. {\bf 19}, 15-18.

\bibitem{Cook1998}
 Cook, R. D. (1998) {\it Regression Graphics: Ideas for
Studying Regressions through Graphics}. New York: Wiley.

\bibitem{}
Cook, R. D. and  Weisberg, S. (1991)  Comment on ``Sliced inverse regression for dimension reduction''. \textit{ Journal of American Statistical Association}, {\bf 86}, 28--33.

%\bibitem{}
% Cortes, C. and Vapnik, V. (1995) Support vector networks. {\it Machine Learning}, {\bf 20}, 1--25.

%\bibitem{}
%Dette, H. and Volgushev, S. (2008). Non-crossing non-parametric estimates of quantile curves. {\it Journal of Royal Statistical Society, Series B}, {\bf 70}, 609--627.

%\bibitem{}
%Efron, B. (1991) Regression percentiles using asymmetric squared error loss. {\it Statistica Sinica}, {\bf 1}, 93--125.

%\bibitem{}  Fukumizu, K.,  Bach, F. R. and Jordan, M. I. (2004)
%Dimensionality reduction for supervised learning with reproducing kernel Hilbert spaces. {\it Journal of Machine Learning Research}, {\bf 5}, 73--79.


%\bibitem{}
%Gu, Y. and Zou, H. (2016) High-dimensional generalizations of asymmetric least squares regression and their applications. {\it The Annals of Statistics}, {\bf 44}, 2661--2694.

\bibitem{}
Harrison, D. and Rubinfeld, D. L. (1978) Hedonic housing prices and the demand for clean air. {\it Journal of environmental economics and management}, {\bf 5}, 81--102.

%\bibitem{}
%Jones, M. C. (1994) Expectiles and M-quantiles are quantiles. {\it Statistics and Probability Letters}, {\bf 20}, 149--153.
%
%\bibitem{}
%Kim, M. and Lee, S. (2016) Nonlinear expectile regression with application to value-at-risk and expected shortfall estimation. {\it Computational Statistics and Data Analysis}, {\bf 94}, 1--19.


\bibitem{}
Kim, H., Wu, Y. and Shin, S. J. (2019) Quantile-slicing estimation for dimension reduction in regression. {\it Journal of Statistical Planning and Inference}, {\bf 198}, 1--12.

\bibitem{}
Li, B. (2018) {\it Sufficient Dimension Reduction: Methods and Applications with R}. CRC Press.


\bibitem{}
Li, B., Artemiou, A. and Li, L. (2011) Principal support vector machines for linear and nonlinear sufficient dimension reduction. {\it  The Annals of Statistics}, {\bf 39}, 3182--3210.


\bibitem{}
Li, B., Kim, M. K. and Altman, N. (2010) On dimension folding of matrix- or array-valued statistical objects. {\it  The Annals of Statistics}, {\bf 38}, 1094--1121.

\bibitem{}
Li, B. and Song, J. (2017) 
Nonlinear sufficient dimension reduction for functional data.  {\it
The Annals of Statistics}, {\bf 45}, 1059--1095.

\bibitem{}
Li, B. and Wang, S. (2007) On directional regression for dimension reduction. \textit{Journal of American Statistical Association},  {\bf 479}, 997--1008.

%\bibitem{} Li, L. and Li, H. (2004) Dimension reduction methods for microarrays with application to censored survival data. {\it Bioinformatics}, {\bf 20}, 3406--3412.




\bibitem{} Li, K. C. (1991) Sliced inverse regression for dimension
reduction (with discussion). {\it Journal of the American
	Statistical Association}, {\bf 86}, 316--342.

\bibitem{}
Li,  L. (2007) Sparse sufficient dimension reduction. {\it Biometrika}, {\bf 94}, 603--613.


%\bibitem{}
%Liu, Y., Chiaromonte, F. and Li, B. (2017) Structured ordinary least squares: a sufficient dimension reduction approach for regressions with partitioned predictors and heterogeneous units. {\it Biometrics}, {\bf 73}, 529--539.

%\bibitem{}
%Liu, Y. and Wu. Y. (2011).
%Simultaneous multiple non-crossing quantile regression estimation using kernel constraints.  {\it Journal of Nonparametric Statistics}, {\bf 23}, 415--437.

%\bibitem{}
%Newey, Whitney K and McFadden, Daniel (1994) Large sample estimation and hypothesis testing. {\it Handbook of econometrics}, {\bf 4}, 2111--2245.

\bibitem{}
Newey, W. K. and  Powell, J. L. (1987) Asymmetric least squares estimation and testing. {\it Econometrica}, {\bf 55}, 819--847.


%\bibitem{}
%Oh, H.S., Nychka, D., Brown, T. and Charbonneau, P., (2004) Period analysis of variable stars by robust smoothing. {\it Journal of the Royal Statistical Society: Series C (Applied Statistics)}, {\bf 53 },15-30.

%\bibitem{}
%Prendergast, L.A., (2008) Trimming influential observations for improved single-index model estimated sufficient summary plots. {\it Computational Statistics \& Data Analysis}, {\bf 52},5319-5327.
%
%\bibitem{}
%Schwarz, Gideon and others (1978) Estimating the dimension of a model. {\it The Annals of Statistics}, {\bf 6}, 461--464.

\bibitem{}
Shin, S. J. and Artemiou A. (2017) Penalized principal logistic regression for sparse sufficient dimension reduction. {\it Computational Statistics and Data Analysis}, {\bf 111}, 48--58.

\bibitem{}
Shin, S. J., Wu, Y., Zhang, H. and Liu, Y. (2017) Principal weighted support vector machines for sufficient dimension reduction in binary classification. {\it Biometrika}, {\bf 104}, 67--81.

\bibitem{} Sz\'ekely, G. J., Rizzo, M. L. and Bakirov, N. K. (2007) Measuring and testing dependence by correlation of distances.
{\it  The Annals of Statistics}, {\bf 35},  2769--2794.


%\bibitem{} Suzuki, T. and Sugiyama, M. (2013)
%Sufficient dimension reduction via squared-loss mutual information estimation.
%{\it Neural Computation}, {\bf 25}, 725--758.

%\bibitem{}
%Sz{\'e}kely, G{\'a}bor J and Rizzo, Maria L and others (2014) Partial distance correlation with methods for dissimilarities. {\it The Annals of Statistics}, {\bf 42}, 2382--2412.

%\bibitem{}
%Taylor, J. W. (2008) Estimating value at risk and expected shortfall using expectiles. {\it Journal of Financial Econometrics}, {\bf 6}, 231--252.
%
%\bibitem{}
%Takeuchi, Ichiro and Nomura, Kaname and Kanamori, Takafumi (2009) Nonparametric conditional density estimation using piecewise-linear solution path of kernel quantile regression. {\it Neural Computation}, {\bf 21}, 533--559.

\bibitem{}
Van der Vaart, A. W. (2000) {\it Asymptotic Statistics}.  Cambridge University Press.

%\bibitem{}
%Vapnik, N. V. (1998) {\em Statistical Learning Theory}. John Wiley \& Sons, Inc.


\bibitem{} Wang, C.,  Shin, S. J.  and Wu, Y. (2018) Principal quantile regression for sufficient dimension reduction with heteroscedasticity.
{\it Electronic Journal of Statistics}, {\bf 12}, 2114--2140.


\bibitem{} 
Wu, H. M. (2008) Kernel sliced inverse regression with applications to classification. {\it Journal of Computational and Graphical Statistics}, {\bf 17}, 590--610.


%\bibitem{}
%Yang, Y., Zhang, T. and Zou, H. (2018)
%Flexible expectile regression in reproducing kernel Hilbert spaces.
%{\it Technometrics}, {\bf 60}, 26--35.

%\bibitem{}
%Yao, Q. and Tong, H. (1996) Asymmetric least squares regression estimation: a nonparametric approach.  {\it Journal of Nonparametric Statistics}, {\bf 6}, 273--292.

\bibitem{}
Yin, X., Li, B. and Cook, R. D. (2008) Successive direction extraction for estimating the central subspace in a multiple-index regression. {\it Journal of Multivariate Analysis}, {\bf 99}, 1733--1757.

%\bibitem{}
%Xia, Y., Tong, H., Li, W. K. and Zhu, L. X. (2002) An adaptive estimation of dimension reduction space. {\it Journal of the Royal Statistical Society: Series B}, {\bf 64}, 363--410.

\end{thebibliography}
\end{document}